\documentstyle[12pt,amscd]{amsart}
\newtheorem{thm}{Theorem}[section]
\newtheorem{lem}[thm]{Lemma}
\newtheorem{prop}[thm]{Proposition}
\newtheorem{cor}[thm]{Corollary}

\newtheorem{remark}[thm]{Remark}
\theoremstyle{definition}
\newtheorem{definition}[thm]{Definition} 

\def\real{{\Bbb R}}
\def\nat{{\Bbb N}}
\def\A{{\cal A}}
\def\B{{\cal B}}
\def\E{{\cal E}}
\def\K{{\cal K}}
\def\L{{\cal L}}

\def\U{{\cal U}}
\def\S{{\cal S}}
\def\chix{\raise.5ex\hbox{$\chi$}}
\def\varep{\varepsilon}
\def\tagto{\mathop{\leftharpoonup\joinrel\rightarrow}\nolimits}
\newcommand{\conv}{\operatorname{conv}}

\begin{document}
\baselineskip=17pt
\title[Complemented copies of $\ell_1$]
{Complemented copies of $\ell_1$ in spaces of vector measures and applications}
\author{Narcisse Randrianantoanina}
\address{Department of Mathematics, The University of Texas at Austin, 
Austin, TX 78712-1082}
\email{nrandri@@math.utexas.edu}
\maketitle

\begin{abstract}
Let $X$ be a Banach space and $(\Omega,\Sigma)$ be a measure space. 
We provide a characterization of sequences in the space of $X$-valued 
countably additive measures on $(\Omega,\Sigma)$ of bounded variation 
that generates complemented copies of $\ell_1$. 
As application, we prove that if a dual Banach space $E^*$ has 
Pe{\l}czy\'nski's property (V*) then so does the space of $E^*$-valued 
countably additive measures with the variation norm. 
Another application, we show that for a Banach space $X$, the space 
$\ell_\infty (X)$ contains a complemented copy  of $\ell_1$ if 
and only if $X$ contains all $\ell_1^n$ uniformly complemented. 
\end{abstract}

\section{Introduction}

Let $X$ be a Banach space and $(\Omega,\Sigma)$ be a measure space. 
We denote by $M(\Omega,X)$ the space of $X$-valued countably additive 
measures on $(\Omega,\Sigma)$ that are of bounded variation with the 
usual total variation norm. 
If $E$ and $F$ are Banach spaces, we denote by $E\widehat\otimes_\pi F$
(resp. $E \widehat\otimes_\varep F$) 
the projective (resp. injective) tensor product of $E$ and $F$. 
We will say that a sequence $(x_n)_n$ is equivalent to a complemented 
copy of $\ell_1$  in $X$ if $(x_n)_n$ is equivalent to the unit vector basis 
of $\ell_1$ and its closed linear span is complemented in $X$. 

In \cite{RAN3} characterization of sequences that are equivalent to 
complemented copy of $\ell_1$ in Bochner spaces were given. 
In this paper we will extend the characterization of \cite{RAN3} 
for the space $M(\Omega,X)$. 
When  $X$ does not have the Radon-Nikodym property one cannot hope to 
represent a measure by its Bochner derivative, but we can always have 
a weak*-density valued in $X^{**}$. 
That is the approach we will take to characterize sequences  in 
$M(\Omega,X)$ that generates complemented  copy of $\ell_1$. 
This approach was used by Talagrand in \cite{T2} to reduce the study
of spaces of vector measures to that of vector valued function spaces.  
 Our main result can be viewed as a common generalization of Theorem~1
of \cite{RAN3} and Theorem~15 of \cite{T2}. We show that if $(e_n)_n$ denotes
the unit vector basis of $c_0$ and $(m_n)_n$ a sequence in $M(\Omega,X)$ then
there exists a sequence $G_n \in \text{conv}\{m_n,m_{n+1},\dots\}$ such that
$\rho(G_n)(\omega) \otimes e_n)_n$ is either weakly Cauchy or equivalent
to the $\ell_1$ basis in $X^{**} \widehat{\otimes}_\pi c_0$ (here 
$\rho(G_n)(\cdot)$
 denoted a weak*-density of $G_n$ that will be described below).
We then use this result to test whether  a given sequence $(m_n)_n$ 
 in $M(\Omega,X)$
contains a subsequence equivalent to a complemented copy of $\ell_1$ or not.

In Section~2, we set some  preliminary background about the space 
$M(\Omega,X)$. 
In particular, we extend the characterization of the dual of 
$M(\Omega,X)$ given by Talagrand in \cite{T2} to characterization of 
operators from $M(\Omega,X)$ into $\ell_1$. 

Section~3 is devoted to the proof of our main theorem about the 
characterization of sequences in $M(\Omega,X)$ that generates 
complemented copy of $\ell_1$. 

In Section~4 and Section~5, we apply the main theorem to study 
Pe{\l}czy\'nski's property (V*) for $M(\Omega,X)$ and complemented 
copies of $\ell_1$ in $\ell_\infty (X)$ and $L^\infty (\lambda,X)$.
We proved that if $X$ is a dual space then $M(\Omega,X)$ has property (V*)
whenever $X$ does. Using the structure of spaces of vector measures, we were 
able to characterize Banach spaces $X$ so that $\ell_\infty(X)$ contains
complemented copy of $\ell_1$. This problem was raised in \cite{SS7} 
(P. 389) (see also \cite{MEN2} Problem~1). Diaz proved in \cite{DI} that 
if $X$ is a Banach lattice then $\ell_\infty(X)$ contains complemented
copy of $\ell_1$ if and only if $X$ contains all $\ell_1^n$ uniformly
complemented. We obtained that similar characterization holds without the
 lattice assumption. 
 
Unexplained notation and terminology can be found in the books 
\cite{D1} and \cite{WO}. 

\section{Preliminaries} 

Let $(\Omega,\Sigma)$ be a measure space and $X$ a Banach space. 
For $m\in M(\Omega,X)$, we denote by $|m|$ its variation. 
Let $\lambda$ be a probability measure on $\Omega$ with $|m|\le\lambda$ 
and consider $\rho$ a lifting of $L^\infty (\lambda)$ \cite{IT}. 
For each $x^* \in X^*$, the scalar measure $x^*\circ m$ has density 
$(d/d\lambda) (x^*\circ m) \in L^\infty (\lambda)$. 
We define $\rho (m)(\omega)$ to be the element in $X^{**}$ defined by 
$$\rho (m)(\omega) = \rho \left( {d\over d\lambda} (x^*\circ m)\right) 
(\omega)\ .$$ 
It is known that $x^* (m(A)) = \int_A \langle \rho (m)(\omega),x^*
\rangle\,d\lambda (\omega)$ for each measurable subset $A$ of 
$\Omega$ and each $x^*\in X^*$. 
Similarly, it can be shown that 
$$|m|(A) = \int_A \| \rho (m)(\omega)\|\, d\lambda (\omega)
\hbox{ for every } A\in\Sigma\ .$$ 
In the case $X=E^*$ is a dual space, $\rho (m)(\omega)$ will be 
the element of $X$ defined by 
$$\rho (m)(\omega)(y) = \rho \left( {d\over d\lambda} (y\circ m)
\right) (\omega) \hbox{ for every } y\in E\ .$$ 
Let us denote by $\tau$ the product on $(X_1^{**})^\nat$ of the 
weak*-topology. 
It was observed in \cite{T2} that if $(m_n)_n$ is a sequence in 
$M(\Omega,X)$ such that $|m_n|\le \lambda$ for every $n\in\nat$ 
then the map $\theta :\Omega \to (X_1^{**})^\nat$ defined by 
$$\theta (\omega) = (\rho (m_n)(\omega))_{n\ge 1} 
\hbox{ is $\tau$-Borel measurable.}$$

For a Banach space $X$ and a compact Hausdorff space $\Omega$, we
denote by $C(\Omega,X)$ the space of all $X$-valued continuous functions
with domain $\Omega$. It is a well known fact that $C(\Omega,X)^*$ is
  isometrically isomorphic to $M(\Omega,X^*)$. This fact will be used in the
sequel.
 
We also need some duality results between spaces of operators and tensor 
products. Let $X$ and $Y$ be Banach spaces. The space of bounded operators
from $X$ into $Y$ will be denoted by $\L(X,Y)$; $N(X,Y)$ will stand for the 
space of nuclear operators from $X$ into $Y$ and $I(X,Y)$ will denote
the space of integral operators from $X$ into $Y$. We  refer to \cite{DU}
and \cite{WO} for  basic properties of these spaces. The following
 identifications
will be used throughout this paper; for what follows the symbol
$\approx$ means isometrically isomorphic.

\begin{prop}
For Banach spaces $X$ and $Y$,

\begin{itemize}
 \item[(i)] $X^* \widehat\otimes_\pi Y \approx N(X,Y)$;
 \item[(ii)] $(X \widehat\otimes_\pi Y)^* \approx \L(X,Y^*)$;
 \item[(iii)] $(X \widehat\otimes_\varep Y)^* \approx I(X,Y^*)$.
\end{itemize}
\end{prop}

As in \cite{RAN3}, the notation $X_1$ will be used for the closed unit ball
of a given Banach space $X$; and finally, to avoid confusion, we will denote 
by $(e_n)_n$ the unit vector basis of $c_0$ while those of $\ell_1$ will be 
denoted by $(e_n^*)_n$.
 
Following the approach of \cite{T2}, we will study bounded operators 
from $M(\Omega,X)$ into $\ell_1$.

\begin{lem} Let $X$ be a Banach space and
 $\S=\{S^*, \ S\in \L(c_0,X)_1\}$.
Then $\S$ is weak*-dense in $\L(X^*,\ell_1)_1$.
\end{lem}

\begin{pf}
If we denote by $\pi_n:\ell_1 \to \ell_1^n$ the canonical projection
 then for 
every operator $T \in \L(X^*,\ell_1)$, the sequence
 $(\pi_n \circ T)_{n \in \nat}$ converges to $T$ for the strong operator
 topology so $\bigcup_{n\geq 1} \L(X^*,\ell_1^n)_1$ is weak*-dense in
$\L(X^*,\ell_1)_1$.
If we denote by $\S_n =\{S^*,\ S\in \L(\ell_\infty^n,X)_1\}$ then it is
 enough
to show that $\S_n$ is weak*-dense in $\L(X^*,\ell_1^n)_1$. To see this
 notice 
that $K(\ell_\infty^n,X)^{**}=\L(\ell_\infty^n,X^{**})$. Let $T$ be an
 element
of $\L(X^*,\ell_1^n)_1$; $T^{*} \in \L(\ell_\infty^n,X^{**})_1$ so the
 exists
a net $(S_\alpha)_\alpha$
 of elements of $\L(\ell_\infty^n,X)_1$ that converges to $T^*$ for 
the weak*-topology in $\L(\ell_\infty^n,X^{**})_1$ and it is now easy to
 check
that $(S_\alpha^*)_\alpha$ converges to $T$ for the weak*-topology in
$\L(X^*,\ell_1^n)_1$. 
 The lemma is proved.
\end{pf}
 
Let $\B$ be a subset of $\L(M(\Omega,X),\ell_1)_1$ defined as follows: 

$T\in \B$ if and only if there exists a finite partition 
$\Omega_1,\Omega_2,\ldots,\Omega_n$ of $\Omega$ and $t_1,t_2,\ldots,t_n 
\in \L(X,\ell_1)_1$ such that $\forall \ m\in M(\Omega,X)$ 
$$T(m)= \sum_{i=1}^n t_i (m(\Omega_i))\in \ell_1\ .$$ 
Clearly $\B$ is a subset of $\L(M(\Omega,X),\ell_1)_1$ and is convex. 

\begin{lem} 
The set $\B$ is weak*-dense in $\L(M(\Omega,X),\ell_1)_1$. 
\end{lem} 
\begin{pf}
To see the lemma, note that from \cite{T2} p.719, if $E$ is a Banach 
space and $\A$ is the subset of $M(\Omega,E)^*$ defined by 
$\varphi\in \A$ if and only if there exist a finite measurable 
partition $\Omega_1,\ldots,\Omega_n$ of $\Omega$ and $x_1^*,\ldots,
x_n^*$ in $E_1^*$ such that $\varphi(m) = \sum_{i=1}^n x_i^* 
(m(\Omega_i))$ then $\A$ is weak*-dense in the unit ball of 
$M(\Omega,E)_1^*$. As above, since
 $\bigcup_{n\geq 1} \L(M(\Omega,X),\ell_1^n)_1$ is dense in 
$\L(M(\Omega,X),\ell_1)_1$ for the strong operator topology, it is enough to 
check that for each $n \geq 1$, $\L(M(\Omega,X),\ell_1^n)_1 \subset
\overline{\B}^{\text{weak*}}$.
 Notice also that if $T\in \L(M(\Omega,X),\ell_1)_1$ and there exist
 measurable partition $\Omega_1, \Omega_2, \dots,\Omega_k$ of $\Omega$ and
$t_1, t_2, \dots, t_k$ in $\L(X,\ell_1)$ such that 
$T(m)= \sum_{i \leq k} t_i \left(m(\Omega_i)\right)$ for all
 $m \in M(\Omega,X)$, then $\| t_i\|\leq 1 $ for all $i \leq k$. 

For each $n \in \nat$, we have
 $\L(M(\Omega,X),\ell_1^n) \approx
 (M(\Omega,X)\widehat\otimes_\pi \ell_\infty^n)^*$. We claim that 
the space $M(\Omega,X) \widehat\otimes_\pi \ell_\infty^n$ is isomorphic to
$M(\Omega,X \widehat\otimes_\pi \ell_\infty^n)$. To see this notice that
 every element $u$ of $M(\Omega,X)\widehat\otimes_\pi \ell_\infty^n$
 has a unique
 representation $u=\sum_{i=1}^n m_i \otimes e_i$. Consider the following map
$J: M(\Omega,X)\widehat\otimes_\pi \ell_\infty^n 
 \to M(\Omega,X \widehat\otimes_\pi \ell_\infty^n)$ defined by
$J(\sum_{i=1}^n m_i \otimes e_i)= \sum_{i=1}^n m_i(\cdot)\otimes e_i$.
 The operator
$J$ is well defined, $\|J\| \leq 1$
 and it is immediate that $J$ is one to one and onto. So
$J^*: M(\Omega, X\widehat\otimes_\pi \ell_\infty^n)^* \to 
\L(M(\Omega,X),\ell_1^n)$ is also one to one and onto. There exists a constant
$C$ such that
 $C.J^*\left((M(\Omega,X)\widehat\otimes_\pi\ell_\infty^n)_1^*\right)$
 contains 
$\L(M(\Omega,X),\ell_1^n)_1$. Define $\A$ as above
 ($E=X\widehat\otimes_\pi \ell_\infty^n$). Since $\A$ is weak* dense in 
$M(\Omega,X\widehat\otimes_\pi \ell_\infty^n)_1^*$, the set
$C.J^*(\A) \cap \L(M(\Omega,X),\ell_1^n)_1$ is weak* dense in the space
$\L(M(\Omega,X),\ell_1^n)_1$. We conlude the proof of the lemma by
 noticing
that $C. J^*(\A) \cap \L(M(\Omega,X),\ell_1^n)_1$ is contained in $\B$.
\quad 
The lemma is proved. 
\end{pf}

Let $\lambda$ be a probability measure on $\Omega$ with $|m|\le\lambda$. 
We define a map $Z(\cdot,m)$ in $\B$ as follows: 
$$Z(\cdot,m) :\B \to L^\infty (\lambda,\ell_1)$$ 
with if $\varphi \in \B$ is defined by the finite measurable partition 
$(\Omega_1,\ldots,\Omega_n)$ of $\Omega$ and $t_1,\ldots,t_n \in \L 
(X,\ell_1)_1$, 
$$Z(\varphi,m) = \sum_{i\le n} {d(t_i\circ m)\over d\lambda} 
\chix_{\Omega_i}\ .$$ 
Notice that since $\ell_1$ has the RNP, the measure $t_i\circ m$ has Bochner
density  with respect to $\lambda$, hence 
 $d(t_i\circ m)/d\lambda$ 
exists and belongs to $L^\infty (\lambda,\ell_1)$. 
Also for every measurable subset $C$ of $\Omega$, 
\begin{equation}
\hbox{Bochner-}\int_C Z(\varphi,m)\,d\lambda = 
\varphi (m_C)
\end{equation}
where $m_C (A) = m(A\cap C)$ for every $A\in\Sigma$. 
In particular 
$$\hbox{Bochner-} \int_\Omega Z(\varphi,m)\, d\lambda =
\varphi (m)\ .$$ 
Note that $L^\infty (\lambda,\ell_1)$ is the dual of $L^1(\lambda,c_0)$. 
Let $(\varphi_\alpha)_\alpha \in \B$ such that for every $\E\in 
M(\Omega,X)\widehat\otimes_\pi c_0$, 
$\lim_\alpha \langle \varphi_\alpha,\E\rangle$ exists. 
For every $C\in\Sigma$ and $u\in c_0$, 
\begin{align*} 
\int_C \langle Z(\varphi_\alpha,m),u\rangle\,d\lambda 
& = \langle \varphi_\alpha (m_C),u\rangle\cr 
& = \langle \varphi_\alpha ,m_C\otimes u\rangle\cr 
& = \langle Z(\varphi_\alpha,m),u\chix_C\rangle 
\end{align*} 
so for any simple function $f$ in $L^1(\lambda,c_0)$, 
$\lim_\alpha \langle Z(\varphi_\alpha,m),f\rangle$ exists and 
since simple functions are dense in $L^1(\lambda,c_0)$, 
we conclude that the map $Z(\cdot,m)$ is weak*- uniformly 
continuous.  
Hence it has a continuous extension still denoted by 
$\varphi\mapsto Z(\varphi,m)$ from $\L(M(\Omega,X),\ell_1)_1$ 
into $L^\infty (\lambda,\ell_1)$. 
Equation~(1) is still valid for every 
$\varphi \in \L(M(\Omega,X), \ell_1)_1$. 

\begin{prop}
Let $\varphi\in \L(M(\Omega,X),\ell_1)_1$ and $(m_n)_n$ 
be a sequence in $M(\Omega,X)$ with $|m_n| \le\lambda$ for 
every $n\in\nat$. 
There exists a countable subset $D$ of the unit ball of 
$\L(X^{**},\ell_1)$ and a map $\omega\mapsto T(\omega) \in 
\overline{D}^{\sigma (\L(X^{**},\ell_1),X^*\widehat\otimes_\pi 
c_0)}$ so that for every $n\in\nat$, 
$$Z(\varphi,m_n) (\omega) = T(\omega) (\rho (m_n)(\omega))
\hbox{\rm  a.e.}$$ 
\end{prop}

\begin{pf}
Here we adopt the methods used in \cite{T2} and \cite{RS2} (Lemma~1). 
Recall that  $(e_j)_j$ denotes  the unit vector basis of $c_0$. 
For each $j\in\nat$ and $m\in M(\Omega,X)$, $|m|\le \lambda$, 
we define 
$$Z_j (\cdot,m): \L(M(\Omega,X),\ell_1)_1 \to L^1(\lambda)$$ 
by $Z_j (\varphi,m)(\omega) = \langle Z(\varphi ,m)(\omega),e_j\rangle$ 
a.e.
 
The map $Z_j(\varphi,m)$ takes every $\varphi \in \L(M(\Omega,X), 
\ell_1)_1$ into $L^\infty (\lambda)$ and since in the unit ball of 
$L^\infty (\lambda)$, the topologies $\sigma (L^\infty,L^1)$ 
and $\sigma (L^1,L^\infty)$ coincide, the map $Z_j(\cdot,m)$ is 
weak* to weakly continuous. 

For each $p\in \nat$, define 
\begin{align*}
Z^{(p)} :\L(M(\Omega,X),\ell_1)_1 \longrightarrow &\ (L^1(\lambda))^{2p}\cr 
\varphi \tagto &\  (Z_j (\varphi,m_n))_{\scriptstyle n\le p\atop\scriptstyle 
j\le p} 
\end{align*}
Same as above, $Z^{(p)}$ is weak* to weakly continuous. 
Since $\B$ is weak* dense in $\L(M(\Omega,X),\ell_1)_1$, we can 
choose an operator $\varphi_p\in \B$ such that 
$$\|Z^{(p)} (\varphi_p) - Z^{(p)} (\varphi)\|\le 2^{-p}$$ 
so 
$$\sup \left\{ \| Z_j(\varphi_p,m_n) - Z_j (\varphi,m_n)\|_1,\ 
n\le p,\ j\le p\right\} \le 2^{-p}\ .$$
Hence for each $n\in \nat$ and $j\in\nat$ 
$$\lim_{p\to\infty} \langle Z(\varphi_p,m_n)(\omega) ,e_j\rangle 
= \langle Z(\varphi,m_n)(\omega),e_j\rangle$$ 
a.e. and we deduce that for a.e. $\omega$, 
$$\hbox{weak*-}\lim_{p\to\infty} Z(\varphi_p,m_n)(\omega) 
= Z(\varphi,m_n)(\omega)\ \ \forall\ n\in\nat\ .$$ 
For each $p\in\nat$, let $\varphi_p (m) = \sum_{i\le k_p} T_{i,p} 
(m(\Omega_{i,p}))$ where $T_{i,p} \in \L(X,\ell_1)_1$ and 
$(\Omega_{i,p})_{i\le k_p}$ is a measurable partition of $\Omega$. 
Let $J:c_0\to \ell^\infty$ be the natural inclusion. 
We need the following lemma:

\begin{lem} 
Let $T \in \L(X,\ell_1)$ and $m\in M(\Omega,X)$ with $|m|\le\lambda$.
The Bochner density of the measure $T\circ m\in M(\Omega,\ell_1)$ is 
given by $\omega\mapsto J^*\circ T^{**} (\rho (m)(\omega))$.
\end{lem} 

To see the lemma, notice that for every $u\in c_0$, 
the scalar valued map
 given by $\omega \to \langle J^*\circ T^{**}(\rho(m)(\omega)),u\rangle$
 is measurable. In fact, for each $\omega \in \Omega$,    
$\langle J^* \circ T^{**} (\rho (m)(\omega)),u\rangle = 
\langle \rho (m)(\omega), T^* Ju\rangle$ and since 
$T^* (Ju)\in X^*$, the function $J^*\circ T^{**} (\rho (m)(\cdot))$ 
is weak*-scalarly measurable and since $\ell_1$ is a separable dual, 
$J^* \circ T^{**} (\rho (m)(\cdot))$ is norm-measurable. 
Moreover for every $A\in\Sigma$ and $u\in c_0$, 
\begin{align*}
\Big\langle \int_A J^*\circ T^* (\rho (m)(\omega))\,d\lambda(\omega),u
\Big\rangle 
& = \int_A \langle \rho(m)(\omega),T^*Ju\rangle \,d\lambda(\omega)\cr 
& = \langle m(A),T^*Ju\rangle 
=  \langle T\circ m(A),Ju\rangle.
\end{align*}
The lemma is proved. 

To complete the proof of the proposition, set 
$D=\{ J^* \circ T_{i,p}^{**},$ $i\le k_p$ and $p\ge1\}$. 
The set $D$ is a countable subset of $\L(X^{**},\ell_1)_1$. 
Now consider a free  ultrafilter $\U$ of $\nat$. For each 
$\omega\in\Omega$, let $T(\omega)\in \overline{D}^{\sigma (\L(X^{**}, 
\ell_1),X^{**}\widehat\otimes_\pi c_0)}$ be the 
weak*-limit along $\U$ of the sequence 
$(J^*\circ T_{i(p,\omega),p}^{**})_{p\in\nat}$ where $i(p,\omega)$ 
is the unique $i\le k_p$ so that $\omega\in\Omega_{i,p}$. 
We now have for a.e. $\omega \in \Omega$:
\begin{align*}
Z(\varphi_p,m_n)(\omega) 
& = \ \sum_{i\le k_p} 
{d(T_{i,p}\circ m_n)(\omega) \over d\lambda} 
\ \chix_{\Omega_{i,p}} (\omega) \cr 
& = \ \sum_{i\le k_p} J^*\circ T_{i,p}^{**} (\rho (m_n)(\omega))
\chix_{\Omega_{i,p}} (\omega)\cr 
& = \ \langle J^* \circ T_{i(p,\omega)}^{**} ,\rho (m_n)(\omega)
\rangle 
\end{align*} 
and hence by taking the weak*-limit on $p$, we get 
$$Z(\varphi, m_n) (\omega) = T(\omega) (\rho (m_n)(\omega))\ .$$ 
This completes the proof.
\end{pf}

\begin{remark}
If $X=Y^*$ is a dual space then $\rho (m_n)(\omega)\in X$ for every 
$\omega$. The set $D$ in the above proposition 
can be chosen as a subset of $\L(X,\ell_1)$ 
and the map $\omega\mapsto T(\omega) \in 
\overline{D}^{\sigma (\L(X,\ell_1),X\widehat\otimes_\pi c_0)}$ 
takes its values in $\L(X,\ell_1)_1$. 
\end{remark} 

\section{Complemented copies of $\ell_1$ in $M(\Omega,X)$} 

In this section we will provide a characterization of some sequences 
in $M(\Omega,X)$ that generate complemented copies of $\ell_1$ generalizing 
the characterization given in \cite{RAN3} for Bochner spaces. 
Specifically we will study, for a sequence $(m_n)_n$ in $M(\Omega,X)$, 
the weak convergence of the sequence $(m_n\otimes e_n)_{n\ge1}$ in 
$M(\Omega,X)\widehat\otimes_\pi c_0$. 
For the rest of this section, $(\Omega,\Sigma)$ is a measure space, 
$\lambda$ a probability measure on $(\Omega,\Sigma)$ and $\rho$ a 
lifting of $L^\infty (\lambda)$. 
For simplicity we will denote by 
$$M^\infty(\lambda,X) = \{ m\in M(\Omega,X)\ ,\ |m|\le\lambda\}\ .$$ 
The following result is our main criterion for determining if a 
given sequence in $M(\Omega,X)$ has a subsequence that generates 
a complemented subspace equivalent to $\ell_1$. 

\begin{thm}
Let $(m_n)_n$ be a sequence in $M^\infty(\lambda,X)$. 
Then there exist a sequence $G_n\in \conv \{m_n,m_{n+1},\ldots\}$ and 
two measurable subsets $C$ and $L$ of $\Omega$ with $\lambda (C\cup L)=1$ 
 and such that: 

(a) for $\omega\in C$, $(\rho (G_n)(\omega)\otimes e_n)_{n\ge1}$ 
is weakly Cauchy in $X^{**} \widehat\otimes_\pi c_0$; 

(b) for $\omega\in L$, $(\rho (G_n)(\omega)\otimes e_n)_{n\ge1}$ 
is equivalent to the unit vector basis of $\ell_1$ in $X^{**}
\widehat\otimes_\pi c_0$.
\end{thm} 

\begin{pf}
We will reduce the  proof to the case of Bochner spaces treated in 
\cite{RAN3} based on the approach used by Talagrand in the proof 
of Theorem~13 of \cite{T2}. 
For what follows, we will identify (for a Banach $E$) the space 
$\ell_\infty (E)$ to the space $\L(\ell_1,E)$: 
A sequence $(x_n)_n$ in $\ell_\infty (E)$ is identified to the 
operator $T\in \L(\ell_1,E)$ defined by 
$$T((a_n)_n) = \sum_{n=1}^\infty a_n x_n \hbox{ for every } 
(a_n)_n \in \ell_1\ .$$ 
Let $F$ be a separable Banach space. 
Consider $(U_n)_n$ a countable basis for the weak*-topology 
on $\L(F,\ell_1)_1$. 
As in \cite{RAN3}, we denote by $\K$ the set of all (weak*) compact 
sets of $\L(F,\ell_1)_1$ and we say that a map $\omega\mapsto K(\omega)$ 
from $\Omega$ to $\K$ is measurable if for each $n\in\nat$, the 
set $\{\omega\in\Omega$, $K(\omega)\cap U_n\ne\emptyset\}$ is 
measurable. 

We are now ready to present the proof of theorem. 
Recall that $\rho(m_n)(\omega)\in X^{**}$ for every $n\in\nat$ and 
$\omega\in\Omega$. 
Define for $\omega\in\Omega$, 
$$K_0 (\omega) = \left\{ (t(\rho (m_n)(\omega))_{n\ge1} \ ;\ 
t\in \L (E^{**},\ell_1)_1\right\}\ .$$ 
Clearly $K_0(\omega)$ is a subset of the unit ball of $\ell_\infty (\ell_1)$ 
and by the identification above, we will view $K_0(\omega)$ as a 
subset of $\L(\ell_1,\ell_1)_1$ so $K_0 :\Omega\to \K$ and we claim 
that $K_0(\cdot)$ is measurable. 
To see this for each $k,j\in\nat$ and $\varep>0$, let 
$$U_{k,j} (\varep) = \{ T\in \L(\ell_1,\ell_1)_1\ ,\ 
|\langle Te^*_k ,e_j\rangle | < \varep \}\ .$$ 
We get that 
\begin{align*} 
\{\omega :K_0(\omega) \cap U_{k,j}(\varep) \ne \emptyset\} 
& = \{ \omega :\exists\ T\in K_0(\omega) \cap U_{k,j}(\varep)\}\cr 
& = \{ \omega :\exists\ t\in \L(X^{**},\ell_1)_1\ ;\  
|\langle t(\rho(m_k)(\omega),e_j\rangle |<\varep\} \ .
\end{align*} 
Let $\S= \{S^*, S\in \L(c_0,X^*)_1\} \subseteq \L(X^{**},\ell_1)_1$. 
$\S$ is weak*-dense in $\L (X^{**},\ell_1)_1$ and 
\begin{align*} 
\{\omega :K_0(\omega) \cap U_{k,j}(\varep)\ne\emptyset\}
&= \{\omega : \exists\ t\in \S\ ,\ |\langle t(\rho (m_k)(\omega)),e_j
\rangle |<\varep\}\cr 
&= \{\omega :(\rho (m_n)(\omega))_{n\ge1} \in V_{k,j}(\varep)\}
\end{align*} 
where the set $V_{k,j}(\varep)$ is defined by 
\begin{align*}
V_{k,j}(\varep) & = \{ (\alpha_n)_n \in X_1^{**\nat} ;\exists\ t\in \S\ ,\ 
|\langle t(\alpha_k),e_j\rangle|<\varep\}\cr 
&= \{ (\alpha_n)_n \in X_1^{**\nat}\ ;\ S\in \L(c_0,X^*)_1\ ,\ 
|\langle \alpha_k,Se_j\rangle |<\varep\} 
\end{align*} 
so $V_{k,j}(\varep)$ is an open subset of $X_1^{**\nat}$ for the topology 
$\tau$ (product of the weak*-topology of $X_1^{**}$). 
Now since $\omega\mapsto (\rho (m_n)(\omega))_{n\ge1}$ is 
$\tau$-Borel measurable, the set $\{\omega :K_0(\omega)\cap 
U_{k,j}(\varep) \ne \emptyset\}$ is  measurable and one can conclude that 
$K_0(\cdot)$ is measurable using the fact that 
$(U_{k,j}(\varep))_{k,j,\varep}$ form a subbasis for the  
weak*-topology of $\L(\ell_1,\ell_1)_1$. 

Let $f_n(\omega) = e_n^*$ $\ \forall\ \omega\in \Omega$. 
The sequence $(f_n(\cdot))_n \in L^\infty (\lambda,\ell_1)$. 
We can apply to the sequence $(f_n)_n$ the construction used in the 
proof of Theorem~1 of \cite{RAN3} starting from $h_n^0 =f_n$ and 
$K_0(\omega)$ as above. 
The construction  yields measurable subsets $C$ and $L$ with 
$\lambda (C\cup L)=1$ and a sequence $g_n\in \conv \{f_n,f_{n+1},
\ldots\}$ such that for $\omega\in C$, $\langle t,g_n(\omega)
\otimes e_n\rangle$ converges for every $t\in K_0(\omega)$ 
and for $\omega\in L$, there exists $a,b\in\real$, $a<b$ such that 
for any finite sequence of zeroes and ones $\sigma$, there exists 
$t\in K_0(\omega)$ such that 
\begin{align*}
&\sigma_n = 1 \Rightarrow \langle t,g_n(\omega)\otimes e_n\rangle\ge b\cr 
&\sigma_n = 0 \Rightarrow \langle t,g_n(\omega)\otimes e_n\rangle\le a\ .
\end{align*} 
If $g_n= \sum_{i=p_n}^{q_n} a_i f_i$ with $p_1\le q_1<p_2 \le q_2\cdots$ 
and $\sum_{i=p_n}^{q_n} a_i=1$, the sequence 
$G_n = \sum_{i=p_n}^{q_n} a_i m_i$ satisfies the conclusion of the 
theorem. 
\end{pf}

\begin{thm} 
Let $(m_n)_n$ be a sequence in $M^\infty(\lambda,X)$. 

(a) If for a.e. $\omega\in\Omega$, the sequence $(\rho (m_n)(\omega) 
\otimes e_n)_n$ is weakly Cauchy in $X^{**} \widehat\otimes_\pi c_0$ 
then $(m_n\otimes e_n)_n$ is weakly Cauchy in 
$M(\Omega,X)\widehat\otimes_\pi c_0$; 

(b) If for a.e. $\omega\in\Omega$, the sequence $(\rho(m_n)(\omega) 
\otimes e_n)_n$ is weakly null in $X^{**}\widehat\otimes_\pi c_0$ 
then $(m_n\otimes e_n)_{n\ge1}$ is weakly null in $M(\Omega,X)
\widehat\otimes_\pi c_0$. 

(c) If there is $L\subset \Omega$, $\lambda(L)>0$ such that for 
$\omega\in L$, there exists $k\in\nat$ such that the sequence 
$(\rho (m_n)(\omega)\otimes e_n)_{n\ge k}$ is equivalent to the 
$\ell_1$-basis in $X^{**}\widehat\otimes_\pi c_0$ then there 
exists $k\in\nat$ such that the sequence $(m_n\otimes e_n)_{n\ge k}$ 
is equivalent to the $\ell_1$-basis in $M(\Omega,X) 
\widehat\otimes_\pi c_0$.
\end{thm} 

\begin{pf}
To prove (a), let $T\in \L(M(\Omega,X)\ell^1)_1$. 
By Proposition~2.4,  there exists a map $\omega\mapsto T(\omega) 
(\Omega\mapsto \L(X^{**},\ell_1)_1)$ such that 
$$Z(T,m_n) (\omega) = T(\omega)(\rho (m_n)(\omega))\ \hbox{ a.e. }\ 
\forall\ n\in\nat\ .$$ 
If for a.e. $\omega$, $(\rho (m_n)(\omega)\otimes e_n)_n$ is weakly 
Cauchy in $X^{**} \widehat\otimes_\pi c_0$, then for a.e. $\omega$, 
$$\lim_{n\to\infty} \langle Z(T,m_n)(\omega),e_n\rangle 
= \lim_{n\to\infty} \langle T(\omega)(\rho (m_n)(\omega)),e_n\rangle$$ 
exists for a.e. $\omega\in\Omega$ and therefore 
$$\lim_{n\to\infty} \langle T,m_n\otimes e_n\rangle 
= \lim_{n\to\infty} \int_\Omega \langle Z(T,m_n)(\omega),e_n\rangle
\, d\lambda(\omega)$$ 
exists. 
Hence $(m_n\otimes e_n)_n$ is weakly Cauchy in $M(\Omega,X)\widehat 
\otimes_\pi c_0$. 
The proof of (b) follows by the same argument. 

Let us prove (c). 
Here we will adopt the proof of Theorem~15(c) of \cite{T2}. 
For $\omega\in L$ and $k\in\nat$, let $\alpha (k,\omega)$ be the 
best constant $\alpha$ such that $(\rho (m_n)(\omega) 
\otimes e_n)_{n\ge k}$ is $\alpha$-equivalent to $\ell_1$. 
We have 
$$\alpha (k,\omega) = \hbox{Inf }\biggl\{ \Big\| \sum_{n\ge k} 
a_n \rho (m_n)(\omega) \otimes e_n\Big\| / 
\sum_{n\ge k} |a_n| \biggr\}$$ 
where the infimum is taken over all finite sequence of rationals. 
The map $\alpha (k,\cdot)$ is measurable. 
Since $\lim_{k\to\infty} \alpha (k,\omega)>0$ for all $\omega\in L$, 
there exists $L'\subset L$, $\lambda (L')>0$ and $k\in\nat$ such 
that $(\rho (m_n)(\omega)\otimes e_n)_{n\ge k}$ is $\alpha$-equivalent 
to $\ell_1$ for each $\omega\in L'$. 

Let $P$ and $Q$ be 2-disjoint finite subsets of $[k,\infty)$. 
For each $\omega\in L'$, there exists $T(\omega) \in \L(X^{**},\ell_1)_1$ 
such that 
$$\langle T(\omega) (\rho (m_q)(\omega)),e_q\rangle > \alpha/2 
\ \hbox{ for }\ q\in P$$ 
and 
$$\langle T(\omega) (\rho (m_q)(\omega)),e_q\rangle < -\alpha/2 
\ \hbox{ for }\ q\in Q\ .$$ 
Since $P\cup Q$ is finite, and the set 
$\{ S^*; S\in \L(c_0,X^*)_1\}$ is weak*-dense in 
$\L(X^{**},\ell_1)_1$, the operator $T(\omega)$ above can be chosen 
in such a way that $T(\omega)=S(\omega)^*$ where 
$S(\omega) \in \L(c_0,X^*)_1$. 
Define 
\begin{align*}
U(\omega) 
& = \{ \omega' \in L',\langle T(\omega)(\rho (m_q)(\omega')),
e_q\rangle > \alpha/2,\ q\in P\cr 
&\hbox{and }\ \langle T(\omega),(\rho(m_q)(\omega')),e_q\rangle 
< - \alpha/2,\ q\in Q\ \}\cr 
& = \{ \omega' \in L',\langle \rho (m_q)(\omega'), S(\omega)
e_q\rangle > \alpha/2,\ q\in P\cr 
&\hbox{and }\ \langle \rho(m_q)(\omega'),S(\omega)e_q\rangle 
< - \alpha/2,\ q\in Q\ \}\ .
\end{align*}
We notice as in \cite{T2} that $U(\omega) \subset\rho (U(\omega))$ 
and since $\omega \in U(\omega)$, 
$\rho (U(\omega))\ne\emptyset$ and hence $\lambda(U(\omega))>0$. 
Let $D$ be a countable subset of $L'$ such that $\lambda (\bigcup_{\omega
\in D} U(\omega))$ is maximal (among the possible choices of $D$). 
Then $\lambda (L'\setminus \bigcup_{\omega\in D} U(\omega))=0$, 
so we can choose $\omega_1,\ldots,\omega_n\in L'$ such that 
$$\lambda \biggl(\bigcup_{i\le n} U(\omega_i)\biggr) \ge \lambda 
(L')/2\ .$$ 
Let $\Omega_i = U(\omega_i)\setminus \bigcup_{\ell<i} U(\omega_\ell)$, 
$\Omega_1,\ldots,\Omega_n$ is a measurable partition of $\bigcup_{i\le n} 
\Omega_i$. We define $T\in \L (M(\Omega,X),\ell_1)$ by 
$$T(m) = \sum_{i\le n} T(\omega_i) (m(\Omega_i))\ \forall\ 
m\in M(\Omega,X)\ .$$ 
The operator $T$ is well defined and $\|T\|\le 1$. 
Moreover for 
\begin{align*} 
q\in P\ \Rightarrow \quad&\langle T(m_q),e_q\rangle  > {\alpha\over  4} 
\lambda (L')\cr 
q\in Q\  \Rightarrow \quad&\langle T(m_q),e_q\rangle > -{\alpha\over 4} 
\lambda (L').
\end{align*}
In fact for every $q\in P$, 
\begin{align*}
\langle T(m_q),e_q\rangle 
&= \sum_{i\le n} \langle T(\omega_i)(m_q(\Omega_i)),e_q\rangle\cr 
&= \sum_{i\le n} \langle \hbox{weak*}-\int_{\Omega_i} \rho (m_q) 
(\omega)\,d\lambda,(\omega),  S(\omega_i)e_q\rangle \cr
&= \sum_{i\le n} \int_{\Omega_i} \langle \rho (m_q)(\omega),
S(\omega_i) e_q\rangle \,d\lambda (\omega)\cr 
&> \sum_{i\le n} \lambda (\Omega_i) {\alpha\over2} > {\alpha\over4} 
\lambda (L')\ .
\end{align*} 
The same estimate is valid for $q\in Q$. 

It follows from Rosenthal's argument in \cite{R1} (see also 
\cite{D1} p.205) that $(m_n\otimes e_n)_{n\ge k}$ is 
${\alpha\over4} \lambda (L')$-equivalent to the $\ell_1$-basis in 
$M(\Omega,X)\widehat \otimes_\pi c_0$. 
The theorem is proved. 
\end{pf}

\begin{remark}
If the Banach space $X$ is a dual space, then Theorem 3.1 and 
Theorem~3.2 are still valid with $X^{**} \widehat \otimes_\pi c_0$ 
replaced by $X\widehat\otimes_\pi c_0$. 
The proofs for this case are just notational adjustments  of the 
proofs given here for Theorem~3.1 and Theorem~3.2 respectively. 
\end{remark} 

\section{Pe{\l}czy\'nski's property (V*) for $M(\Omega,X)$}    

\begin{definition}
Let $E$ be a Banach space. A series $\sum_{n=1}^\infty x_n$ in $E$ is 
said to be a Weakly Unconditional Cauchy (W.U.C.) if for every $x^*$ 
in $E^*$, the series $\sum_{n=1}^\infty |x^* (x_n)|$ is convergent. 
\end{definition} 

There are many criteria for a series to be a W.U.C.\ series (see for 
instance \cite{D1} or \cite{WO}). 
The following definition introduced by Pe{\l}czy\'nski in \cite{PL1} 
isolates the class of spaces that we would like to study in this section. 

\begin{definition} 
A Banach space $E$ is said to have property (V*) if a subset $K$ of $E$ 
is relatively weakly compact whenever $\lim_{n\to\infty} \sup_{x\in K}  
|x(x_n^*)|=0$ for every W.U.C.\ series $\sum_{n=1}^\infty x_n^*$ in $E^*$. 
\end{definition} 

\begin{definition}
A subset $K$ of a Banach space $E$ is called a (V*)-set if for every 
W.U.C.\ series $\sum_{n=1}^\infty x_n^*$ in $E^*$, the following holds: 
$\lim_{n\to\infty} \sup_{x\in K} |x(x_n^*)|=0$.
\end{definition} 
Hence a Banach space $E$ has property (V*) if and only if every 
(V*)-set in $E$ is relatively weakly compact.
 
The notion of (V*)-set was introduced and studied by Bombal in \cite{B1}.
 He proved  ( see \cite{B1} Proposition~1.1) that a subset $K$ of 
 a Banach space
is a (V*)-set if and only if for every operator $T\in \L(E,\ell_1)$, 
$T(K)$ is relatively compact in $\ell_1$.
Using the fact that $\L(E,\ell_1)\approx (E\widehat\otimes_\pi c_0)^*$,  the 
following lemma is immediate

\begin{lem}
A subset $K$ of a Banach space $E$ is a (V*) subset if and only if 
for every sequence $(x_n)_n$ in $K$, the sequence 
$(x_n \otimes e_n)_n$ is weakly null in $E \widehat\otimes_\pi c_0$.
\end{lem}
 
Property (V*) has been considered by several authors. 
We refer to \cite{B1}, \cite{EM4}, \cite{GOS2}, \cite{LEU}, \cite{PF1} and 
\cite{RAN5} for more additional information on property (V*). 
Perhaps the most appealing fact about property (V*) is its connection 
with complemented copy of $\ell_1$. 
>From a result of Emmanuele \cite{EM4} (see also \cite{GOS2}), 
the following simple characterization can be deduced. 

\begin{prop} 
A Banach space $E$ has property (V*) if and only if $E$ is weakly 
sequentially complete and every sequence that is equivalent to the 
unit vector basis of $\ell_1$ in $E$ has a subsequence equivalent 
to a complemented copy of $\ell_1$. 
\end{prop} 

The most notable examples of Banach spaces with property (V*) 
are $L^1$-spaces. 
It was shown in \cite{RAN3} that Bochner spaces $L^1(\lambda,X)$ 
has property (V*) whenever $X$ does. 
In this section, we will study the property (V*) for the space 
$M(\Omega,X)$. 
Since $M(\Omega)$ is an $L^1$-space, it has property (V*) and it is a 
natural question to address if $M(\Omega,X)$ has property (V*) 
whenever $X$ does. 
It was shown in \cite{T2} that if $X$ is a weakly sequentially complete 
dual space, then $M(\Omega,X)$ is weakly sequentially complete. 
A counterexample by Talagrand in \cite{T5} reveals that there exists 
a Banach lattice $X$ with property (V*) with $M(\Omega,X)$ containing $c_0$. 
The following theorem can be viewed as a common  
generalization of Theorem~2 of 
\cite{RAN3} and Theorem~17 of \cite{T2}. 

\begin{thm} 
Let $X=Y^*$ be a dual space. 
The space $M(\Omega,X)$ has property (V*) if and only if $X$ does. 
\end{thm}

\begin{pf} 
We will prove the non-trivial implication. 
Let $K\subset M(\Omega,X)_1$ be a (V*)-set. 
The set $K$ does not contain any sequence that generates a 
complemented copy of $\ell_1$. 

Notice that the space $M(\Omega,X)=M(\Omega,Y^*)$ is isometrically 
isomorphic to $C(\Omega,Y)^*$ and therefore $M(\Omega,Y^*)
\widehat\otimes_\pi c_0$ can be identified to the space 
$N(C(\Omega,Y),c_0)$.

The following lemma will be used. 

\begin{lem} 
Let $Y$ be a Banach space and $(y_n^*)_n$ a bounded sequence in $Y^*$. 
Assume that: 

(i) $(y_n^*\otimes e_n)_{n\ge1}$ is weakly Cauchy in $Y^*\widehat\otimes_\pi 
c_0$; 

(ii) There exists a convex combination of $(y_n^* \otimes e_n)_{n\ge1}$ 
say $(\sum_{i=p_n}^{q_n} \alpha_i y_i^*\otimes e_i)_{n\ge1}$ such that 
$$\lim_{n\to\infty} \Big\|\sum_{i=p_n}^{q_n} \alpha_i y_i^* 
\otimes e_i\Big\|_{I(Y,\ell_\infty)} =0$$ 
then $(y_n^* \otimes e_n)_{n\ge1}$ is weakly null in $Y^*\widehat\otimes_\pi 
c_0$. 
\end{lem}

To see the lemma, let $u$ be an element of $Y^*\widehat\otimes_\pi c_0$; 
$u$ can be viewed as an operator from $Y$ into $c_0$. 
If $u$ is of the form $u=\sum_{i=p_n}^{q_n} \alpha_i y_i^*\otimes e_i$ 
then we claim that 
$$\|u\|_{Y^*\widehat\otimes_\pi c_0} = \|u\|_{I(Y,\ell_\infty)}\ .$$ 
To see this notice that $(Y^* \widehat\otimes_\pi c_0)^* = \L(Y^*,\ell_1)$ 
and $I(Y,\ell_\infty)$ is a dual space with predual 
$Y\widehat\otimes_\varep \ell_1 \approx K_{w^*} (Y^*,\ell_1)$ the space of 
weak* to weakly continuous compact operators from $Y^*$ into $\ell_1$. 

Let $\varep>0$, the space $\S=\{S^*,S\in \L(c_0,Y)_1\}$ is weak*-dense 
in $\L(Y^*,\ell_1)_1$ so it is norming. 
There exits $S\in \L(c_0,Y)$, $\|S\|=1$ such that 
$$\|u\|_{Y^*\widehat\otimes_\pi c_0} - \varep 
\le \langle u,S^*\rangle  
= \sum_{i=p_n}^{q_n} \alpha_i\langle y_i^*,Se_i\rangle\ ;$$ 
define $K:c_0\to Y$ by $Ke_i=Se_i$ for $i\in [p_n,q_n]$ and $Ke_i=0$ 
otherwise. 
 The operator $K$ is clearly bounded with $\|K\| \le \|S\|\le 1$. 
 Moreover the operator $K^*\in K_{w^*} (Y^*,\ell_1)$ and 
$\langle u,S^*\rangle = \langle u,K^*\rangle \le \|u\|_{I(Y,\ell_\infty)}$ 
so $\|u\|_{Y^*\widehat\otimes_\pi c_0} \le \|u\|_{I(Y,\ell_\infty)}+\varep$ 
and since $\varep$ is arbitrary, the claim follows.
 
We finish the proof of the lemma by noticing that the sequence $(y_n^* 
\otimes e_n)_n$ is weakly Cauchy and $0$ is a weak-cluster point of 
$(y_n^*\otimes e_n)_n$ so $(y_n^*\otimes e_n)_n$ is weakly null.
 The lemma is proved.

We are now ready to proceed for the theorem. Without loss of generality
we can assume that there exists a probability measure $\lambda$ on 
$(\Omega,\Sigma)$ such that $|m|\leq \lambda$ for every $m \in K$.  
Let $(m_n)_n$ be a sequence in $K$. 
Since $K$ is a (V*)-set, for any sequence 
$G_n\in \conv \{m_n,m_{n+1},\ldots\}$,
Lemma~4.4 implies that   $(G_n\otimes e_n)_n$ is weakly 
null in $M(\Omega,Y^*)\widehat\otimes_\pi c_0$. 
Applying Theorem~3.1 to  the sequence $(m_n)_{n}$, there exist a sequence 
$G_n\in \conv \{m_n,m_{n+1},\ldots\}$, measurable subsets $C$ and $L$ of 
$\Omega$  
 with $\lambda(C\cup L)=1$  and such that: 

(a) for $\omega \in C$, $(\rho (G_n) (\omega) \otimes e_n)_{n\ge1}$ is 
weakly Cauchy in $Y^* \widehat\otimes_\pi c_0$; 

(b) for $\omega\in L$, $(\rho (G_n)(\omega)\otimes e_n)_{n\ge1}$ is 
equivalent to $\ell_1$ in $Y^*\widehat\otimes_\pi c_0$. 

But since $(G_n\otimes e_n)_n$ is weakly null in $M(\Omega,Y^*) 
\widehat \otimes_\pi c_0$, we conclude from Theorem~3.2(c) that 
$\lambda (L)=0$ so for a.e. $\omega$, the sequence $(\rho (G_n)(\omega) 
\otimes e_n)_{n\ge1}$ is weakly Cauchy in $Y^* \widehat\otimes_\pi c_0$. 

On the other hand, the space $M(\Omega,Y^*)\widehat\otimes_\pi c_0$ can 
be identified to $N(C(\Omega,Y),c_0)$. 
Let $J$ be the canonical injection of $N(C(\Omega,Y),c_0)$ into 
$I(C(\Omega,Y),\ell_\infty)$. 
The space $I(C(\Omega,Y),\ell_\infty)$ is isometrically isomorphic to 
$I(C(\Omega),I(Y,\ell_\infty))$ (see Theorem~3 of \cite{S2}) which is 
also isometrically isomorphic to $M(\Omega,I(Y,\ell_\infty))$. 
Let us denote by $\theta :I(C(\Omega,Y),\ell_\infty)\to M(\Omega,I(Y,
\ell_\infty))$ the isomorphism. 

\begin{lem}
Let $G$ be a measure in $M(\Omega,Y^*)$ and $e\in c_0$. 
If we denote by $\widehat G$ the measure in $M(\Omega,I(Y,\ell_\infty))$ 
given by $\widehat G= \theta\circ J(G\otimes e)$ then 
$$\rho (\widehat G)(\omega)=\rho (G)(\omega)\otimes e\ \  a.e.$$ 
\end{lem} 

To see this let $f\in C(\Omega)$ and $y\in Y$ 
\begin{align*} 
\langle \widehat G,f\rangle (y) &= \langle \widehat G,f\otimes y\rangle\cr 
&= \langle G\otimes e,f\otimes y\rangle \cr
&= \int_\Omega f(\omega)\langle \rho (G)(\omega),y\rangle e\,d\lambda 
(\omega) 
\end{align*}
and by the definition of $\rho (\widehat G)$, 
$$\langle \widehat G,f\otimes y\rangle = \int_\Omega f(\omega)\rho 
(\widehat G)(\omega)(y)\,d\lambda (\omega).$$ 
So we conclude that for each $y\in Y$, 
$$\rho (\widehat G)(\omega)(y) = (\rho (G)(\omega)\otimes e)(y)\hbox{ a.e. }$$
and the lemma follows. 

To conclude the proof of the theorem, recall that $(G_n\otimes e_n)_n$ 
is weakly null in the space $M(\Omega,Y^*)\widehat \otimes_\pi c_0$. 
>From the identification above, if we denote by
 $\widehat G_n = \theta \circ J 
(G_n\otimes e_n)$ then $(\widehat G_n)_n$ is weakly null in 
$M(\Omega,I (Y,\ell_\infty))$. Using Theorem~16 of \cite{T2},
  we can find $p_1<q_1 <p_2<q_2\cdots 
p_n<q_n$ and sequences $(\alpha_{i,k})$, with 
$\sum_{i=p_n}^{q_n} \alpha_{i,n}=1$ for all $n\in \nat$ such that
for a.e. $\omega \in \Omega$, 
$$\lim_{n\to\infty} \Big\| \sum_{i=p_n}^{q_n} \alpha_{i,n} 
\rho (\widehat G_i)(\omega)\Big\|_{I(Y,\ell_\infty)} =0.$$ 
But $\rho (\widehat G_i)(\omega) = \rho (G_i)(\omega)\otimes e_i$ so 
applying Lemma~4.7, $(\rho (G_n)(\omega)\otimes e_n)_{n\ge1}$ is weakly 
null in $Y^* \widehat \otimes_\pi c_0$ for a.e. $\omega$, and we
can deduce from Lemma~4.4 that  
$\{\rho (G_n)(\omega),n\ge1\}$ is a V*-set in $Y^*$ and since $Y^*$ 
has property (V*), $\{\rho (G_n)(\omega),n\ge1\}$ is relatively 
weakly compact in $Y^*$ (for a.e. $\omega$). 
Hence by Theorem~1 of \cite{RS2}, $K$ is relatively weakly compact. 
The proof is complete. 
\end{pf}

\section{Complemented copies of $\ell_1$ in $\ell_\infty (X)$}

In this section, we will apply Section 2 and Section 3 to investigate 
the following question: 
when does $\ell_\infty (X)$ contain a complemented copy of $\ell_1$? 
A closely related question was raised in \cite{SS7} p.389 (see also
\cite{MEN2} Problem~1). 
The expected answer would be that $X$ contains a complemented copy 
of $\ell_1$. 
But this is not the case since the space $X= (\bigoplus_n \ell_1^n)_{c_0}$ 
does not contain any copy of $\ell_1$ but $\ell_\infty (X) \approx 
(\bigoplus_n \ell_1^n)_{\ell_\infty}$ contains a 
complemented copy of $\ell_1$.  
Diaz \cite{DI} has found an interesting result in this direction using local 
theory. 
Following \cite{DI}, we will use the following definition: 

\begin{definition} 
Let $1\le p\le\infty$. 
A Banach space $X$ is said to be an $S_p$-space or that it contains all 
$\ell_p^n$ uniformly complemented or it contains a complemented local copy 
of $\ell_p$ if there is some constant  
$C\ge1$ such that for every $n\in\nat$ 
there exist operator $J_n\in \L(\ell_p^n,X)$ and $P_n\in \L(X,\ell_p^n)$ 
such that $P_nJ_n = id_{\ell_p^n} $ and $\|P_n\|\cdot \|J_n\| \le C$. 
\end{definition} 

We notice that  a Banach space $X$ is an $S_p$-space if and only if 
its dual is an $S_{p'}$-space ($1/p +1/p'=1$). As a result
a Banach space $X$ is an $S_p$-space if and only if $X^{**}$ is an
$S_p$-space.
 
For the following theorem, $(\Omega,\Sigma,\lambda)$ is a finite 
measure space. 

\begin{thm} 
Let $X$ be a Banach space. 
The following statements are equivalent.
\begin{enumerate}
\item $L^\infty (\lambda,X)$ contains a complemented copy of $\ell_1$; 
\item $\ell_\infty (X)$ contains a complemented copy of $\ell_1$; 
\item $c_0(X)$ contains all $\ell_1^n$ uniformly complemented; 
\item $X$ contains all $\ell_1^n$ uniformly complemented.
\end{enumerate}
\end{thm} 

This theorem was proved in \cite{DI} under the assumption that $X$ is 
a Banach lattice. 
It should be noted that in the  proof given in \cite{DI},
 the lattice assumption was 
used only to show that $(2)\Rightarrow (3)$. 
So we will present here the implication $(2)\Rightarrow (3)$ and we refer 
to \cite{DI} Theorem~1  for the other implications. 
We will divide the proof into two parts. 

\begin{prop} 
Let $(\Omega,\Sigma,\lambda)$ be a finite measure space and $E$ be a Banach 
space. Let $(m_n)_n$ be a bounded sequence in
 $M^\infty(\lambda,E)$. If $E^{**}$ 
does not contain any complemented copy of $\ell_1$ then the sequence 
$(m_n)_n$ does not have any subsequence that generates a complemented 
copy of $\ell_1$ in $M(\Omega,E)$.
\end{prop} 

\begin{pf} 
To prove that $(m_n)_n$ does not contain any subsequence equivalent to 
complemented copy of $\ell_1$, it is enough to show that for every 
operator $T\in \L(M(\Omega,E),\ell_1)$, $\{T(m_n),n\ge1\}$ is 
relatively compact in $\ell_1$. 
To see this let $T\in \L(M(\Omega,E),\ell_1)_1$ and $\rho$ be a lifting 
of $L^\infty (\lambda)$. 
By Proposition~2.4, there exists a map $\omega\mapsto T(\omega) 
(\Omega\to \L(E^{**},\ell_1))_1$ such that 
$$Z(T,m_n) (\omega) = T(\omega) (\rho (m_n)(\omega))\ \hbox{ a.e.}$$ 
where $Z(\cdot,m_n)$ is the notation introduced in Section~2. 

Recall that $T(m_n)= \int_\Omega Z(T,m_n)(\omega) \,d\lambda (\omega)$. 
Since $E^{**}$ does not contain any complemented copy of $\ell_1$, the 
operators $T(\omega)$'s are compact for every $\omega\in\Omega$ and 
therefore $\{T(\omega)(\rho (m_n)(\omega)),n\ge1\}$ is compact for 
a.e. $\omega$. 
By \"Ulger's criteria of weak-compactness in Bochner spaces, 
$\{T(\cdot)(\rho (m_n)(\cdot)),n\ge1\}$ is relatively weakly compact 
in $L^1(\lambda,\ell_1)$ and hence there exists subsequence $(m_{n_k})$ 
of $(m_n)_n$ such that 
$$T(m_{n_k})= \int_\Omega
 T(\omega)(\rho (m_{n_k})(\omega))\,d\lambda(\omega)$$ 
converges weakly in $\ell_1$. 
The proof is complete.
\end{pf} 

For the next proposition, we consider $\Omega= [0,1]$, 
$\Sigma$ is the $\sigma$-algebra of the Borel subsets of $[0,1]$ and 
$\lambda$ the Lebesgue measure. 

\begin{prop} 
Let $X$ be a Banach  space and let us assume that there exists a bounded 
sequence $(a_k)_{k\in\nat}$ in $\ell_\infty (X)$ that generates a 
complemented copy of $\ell_1$ in $\ell_\infty (X)$ then there exists a 
sequence $(m_k)_k$ in $M([0,1],c_0(X))$ with $m_k\in
 M^\infty(\lambda,c_0(X))$ 
for every $k\in\nat$ and the sequence
 $(m_k)_k$ is equivalent to a complemented copy 
of $\ell_1$ in $M([0,1],c_0(X))$. 
\end{prop} 

\begin{pf} Let $(a_k)_{k\in\nat}$ be a sequence that
 is equivalent to a complemented copy of 
$\ell_1$ in $\ell_\infty (X)$. There exists a bounded operator $S$ from 
$\ell_\infty (X)$ onto $\ell_1$ with $S(a_k)= e_k^* $ $\forall\ k\in\nat$. 
Let $(r_n)_{n\ge1}$ be the sequence of  Rademacher functions. 
For each $k\in\nat$, we define 
$$m_k(A) = \biggl(\int_A r_n(t) a_k(n)\,dt\biggr)_{n\ge1}$$ 
where $a_k(n)$ denotes the $n^{th}$ projection of $a_k\in\ell_\infty (X)$ 
into $X$. 

For each $n\in\nat$, $m_k(A)(n)=\int_A r_n(t)a_k(n)\,dt$ belongs to $X$ and 
$$\|m_k (A)(n)\| = \Big| \int_A r_n(t)\,dt\Big|\ \|a_k(n)\|$$ 
and since $(r_n)_n$ is weakly null in $L^1[0,1]$ (see for instance 
\cite{D1}), 
$\lim\limits_{n\to \infty} \|m_k(A)(n)\|=0$. So $m_k$ is indeed a measure in 
$M([0,1],c_0(X))$. 
Moreover for every measurable subset $A$ of $[0,1]$,
 $\|m_k(A)\|_{c_0(X)} \le \sup_k \|a_k\|.\lambda(A)$ so the 
sequence $(m_k)_k$ is bounded and is in $M^\infty(\lambda,c_0(X))$.
 
For any measure $m\in M([0,1],c_0(X))$, we will denote by $m^{(n)}$ the 
measure in $M([0,1],X)$ given by $m^{(n)}(A) = m(A)(n)$ (the $n^{th}$ 
projection of $m(A)$ onto $X$). 
The measure $m^{(n)}$ is well defined for every $n\in\nat$. 

If $|m|\ll \lambda$, then $|m^{(n)}|\ll \lambda$ for every $n\in\nat$ and 
we will denote by $\rho (m^{(n)}):[0,1]\to X^{**}$ a weak*-density of 
$m^{(n)}$ with respect to the Lebesgue measure $\lambda$. 
In this case 
$$m^{(n)} (A) = \hbox{weak*}-\int_A \rho (m^{(n)})(t)\,dt\in X$$ 
for every $A\in\Sigma$. 
Since $r_n$ is a simple function on $[0,1]$, we conclude that 
$$\hbox{weak*}-\int_A \rho (m^{(n)})(t)r_n(t)\,dt\in X$$ 
for every $n\in\nat$ and $A\in \Sigma$ so the measure $[m]:\Sigma\to 
\ell_\infty (X)$ defined by 
$$[m](A) = \biggl( \hbox{weak*}-\int_A\rho (m^{(n)})(t)r_n(t)\,dt
\bigg)_{n\ge1}$$ 
is well defined and satisfies 
$$\|[m](A)\|_{\ell_\infty(X)} \le \|m(A)\|_{c_0(X)}\ .$$ 
Define $T:M_\lambda ([0,1],c_0(X))\to M([0,1],\ell^\infty (X))$ by 
$$T(m) = [m]\ \ \  \forall\ m\in M_\lambda([0,1],c_0(X)) \ .$$
(Here $M_\lambda ([0,1],c_0(X))$ denotes the subset of $M([0,1],c_0(X))$ 
of all measures that are absolutely continuous with respect to $\lambda$.) 
The operator $T$ is clearly linear and $\|T\|\le1$. 
For the operator $S:\ell_\infty (X)\to \ell_1$ defined above, we consider 
$\widehat S:M([0,1],\ell_\infty(X))\to \ell_1$ by 
$\widehat S(m) = S(m([0,1]))$.

\noindent {\bf Claim:}
For all $k\in\nat$, $\widehat S\circ T(m_k)=e_k^*$.

To see the claim, notice that since 
\begin{align*}
m_k(A) & = \biggl( \int_A r_n(t) a_k(n)\,dt\biggr)_{n\ge1}\ ,\cr 
m_k^{(n)} (A) & = \hbox{\rm Bochner }\int_A r_n(t)a_k(n)\,dt\ .
\end{align*} 
So
\begin{align*} 
T(m_k) (A) = [m_k](A) & = \bigl( \lambda (A)a_k(n)\bigr)_{n\ge1}\cr 
&= \lambda (A) a_k\ .
\end{align*}
Hence $\widehat S\circ T(m_k) = S(a_k)=e_k^*$\  $\forall k\in\nat$. 
This shows that $(m_k)_{k\ge1}$ contains a subsequence equivalent to the 
$\ell_1$-basis in $M_\lambda([0,1],c_0(X))$ and since $M_\lambda ([0,1],
c_0(X))$ is complemented in $M([0,1],c_0(X))$ (see Theorem~I.9 of 
\cite{DU}), the proof of the proposition is complete.
\end{pf} 

To prove the theorem, assume that
 $\ell_\infty(X)$ contains a complemented copy of 
$\ell_1$, then from Proposition~5.4 there exists a bounded sequence 
$(m_k)_k$ in $M^\infty(\lambda,c_0(X))$ that is equivalent to a
 complemented 
copy of $\ell_1$ in $M([0,1],c_0(X))$ and this implies from 
Proposition~5.3 that $c_0(X)^{**}$ contains a complemented copy of
 $\ell_1$ 
so $c_0(X)^{**}$ is an $S_1$-space and this is equivalent to $c_0(X)$ 
is an $S_1$-space. This  completes the proof of the theorem. \quad \qed

\begin{cor}
If $X$ is a super-reflexive space then $\ell_\infty(X)$ does not
contain any complemented  copy of $\ell_1$.
\end{cor}

As it was noticed  in \cite{DI}, the space
$X=\left( \oplus_n \ell_1^n \right)_{\ell_2}$ is a reflexive space
but $L^\infty([0,1],X)$ and $\ell_\infty(X)$ contain complemented copies
of $\ell_1$.

\vskip .5truein
\noindent {\bf Acknowledgements.} 
Most of the work  in this paper was done when the author was 
participating in the workshop on linear analysis and probability 
(Summer 1995) at the Texas A\&M University. The author  would like to 
express his gratitude to the Department of Mathematics of the Texas 
A\&M University for its hospitality and  financial support.


\end{document}